\documentclass{amsart}

\usepackage{amscd,amssymb,amsmath}
\usepackage{graphics}

\newtheorem{theorem}{Theorem}[section]

\newtheorem{lemma}[theorem]{Lemma}
\newtheorem{proposition}[theorem]{Proposition}
\newtheorem{corollary}[theorem]{Corollary}

\theoremstyle{definition}

\newtheorem{remark}[theorem]{Remark}

\theoremstyle{remark}

\renewcommand{\labelenumi}{(\roman{enumi})}

\DeclareFontFamily{U}{wncy}{}
\DeclareFontShape{U}{wncy}{m}{n}{<->wncyr10}{}
\DeclareSymbolFont{mcy}{U}{wncy}{m}{n}
\DeclareMathSymbol{\Sh}{\mathord}{mcy}{"58}

\newcommand\mylabel[1]{\label{#1}}

\newcommand{\ZZ}{\mathbb{Z}}

\newcommand{\RR}{\mathbb{R}}
\newcommand{\CC}{\mathbb{C}}

\newcommand  {\shC}     {\mathcal{C}}

\newcommand  {\shF}     {\mathcal{F}}

\newcommand  {\shH}     {\mathcal{H}}

\newcommand  {\shS}     {\mathcal{S}}

\newcommand  {\foU}     {\mathfrak{U}}

\newcommand  {\crs}     {{\rm crs}}

\newcommand  {\dirlim}  {\varinjlim}

\newcommand  {\lra}     {\longrightarrow}

\newcommand  {\quadand} {\quad\text{and}\quad}

\newcommand  {\ra}      {\rightarrow}

\def\mydate{\number\day\space\ifcase\month \or January\or February\or March\or 
April\or May\or June\or July\or
August\or September\or October\or November\or December\fi \space\number\year}

\DeclareFontFamily{U}{wncy}{}
\DeclareFontShape{U}{wncy}{m}{n}{<->wncyr10}{}
\DeclareSymbolFont{mcy}{U}{wncy}{m}{n}
\DeclareMathSymbol{\Sh}{\mathord}{mcy}{"58}

\begin{document}

\title[Pathologies in cohomology]
      {Pathologies in cohomology of non-paracompact Hausdorff spaces}

\author[Stefan Schr\"oer]{Stefan Schr\"oer}
\address{Mathematisches Institut, Heinrich-Heine-Universit\"at,
40204 D\"usseldorf, Germany}
\curraddr{}
\email{schroeer@math.uni-duesseldorf.de}


\dedicatory{Revised version, 21 June 2013}

\begin{abstract}
We construct  a non-paracompact Hausdorff space for which
\v{C}ech cohomology does not coincide with sheaf cohomology. 
Moreover, the sheaf of continuous real-valued functions is neither soft
nor acyclic, and our space admits non-numerable principal bundles.
\end{abstract}

\maketitle
\renewcommand{\labelenumi}{(\roman{enumi})}

\section*{Introduction}
Recall that a topological space $X$ is called \emph{paracompact} if 
it is Hausdorff, and each open covering admits a refinement that is locally finite.
This notion was introduced by Dieudonn\'e \cite{Dieudonne 1944} as early as 1944 and has turned out
to be extremely useful in general topology and sheaf theory.
For example, Godement showed that   \emph{\v{C}ech cohomology}
coincides with \emph{sheaf cohomology} on paracompact spaces (\cite{Godement 1964}, Theorem 5.10.1).
For general spaces, all that can be said is that there is a spectral sequence
$$
\check{H}^p(X,\shH^q(\shF))\Longrightarrow H^{p+q}(X,\shF)
$$
computing  the ``true'' sheaf cohomology from  the \v{C}ech cohomology of the presheaves of 
sheaf cohomology (loc.\ cit., Theorem 5.9.1).
Grothendieck observed that for many irreducible spaces, for example  $X=\CC^2$ with
the Zariski topology, this spectral sequence does not degenerate for suitable  $\shF$, such that 
\v{C}ech cohomology does not coincide with sheaf cohomology (\cite{Grothendieck 1957}, page 178).
On the other hand, Artin \cite{Artin 1971} established that for ``most''  separated schemes, \v{C}ech cohomology
agrees with sheaf cohomology when computed in the \'etale topology.

Although the known counterexamples are very common in the realm of algebraic geometry,
they are perhaps not so natural from the standpoint of algebraic or general topology, since the spaces are not Hausdorff.
In my opinion, it would be desirable to have further counterexamples satisfying the  Hausdorff axiom,
the more so in light of Artin's result.

The goal of this note is to provide such a space. 
The construction   roughly goes as follows:
We start with   an infinite wedge sum $X=\bigvee_{i=1}^{\infty} D^2$ of closed 2-disks, and replace the CW-topology 
at the intersection of the 2-disks by some coarser topology. This topology is choose fine enough to keep the space
Hausdorff, yet coarse enough so that a variant of Grothendieck's argument holds true.

It turns out that our space has other pathological features as well: The sheaf of continuous real-valued functions is neither soft
nor acyclic. Although the space is contractible, it carries nontrivial principal $S^1$-bundles. These are necessarily
non-numerable, whence do not come from the universal bundle.

\section{The construction}

We start by constructing an infinite 2-dimensional CW-complex $X$.
Its 0-skeleton is a sequence $e_n^0$, $n\geq 0$ of 0-cells.
The first 0-cell $x=e_0^0$ will play a special role throughout, and we shall call
it the \emph{origin}.
To form the 1-skeleton $X^1$, we connect the origin to each $e_n^0$, $n\geq 1$ with two 1-cells called $e_{\pm n}^1$.
To complete the construction, we  choose homeomorphisms
$$
\varphi_n:S^1\lra e_0^0\cup e_n^0\cup e_{n}^1\cup e_{-n}^1\subset X^1,
$$
and use these as attaching maps for the 2-cells $e_n^2$, $n\geq 1$.  This gives an infinite 2-dimensional CW-complex $X$,
which one may visualize as follows.

\vspace{2em}
\centerline{\includegraphics{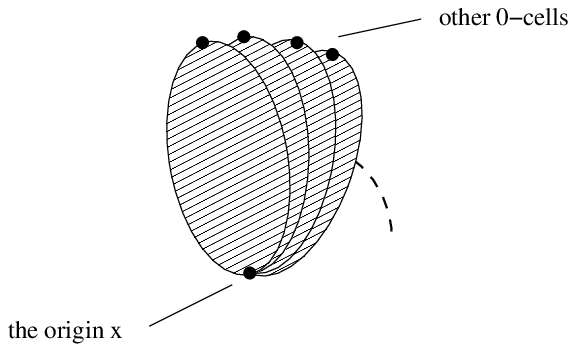}} 
\vspace{1em}
\centerline{Figure \stepcounter{figure}\arabic{figure}: The CW-complex $X$}
\vspace{1em}

Being a CW-complex, the space $X$ is paracompact \cite{Miyazaki 1952}. With our goal in mind we  now replace the CW-topology
by  some coarser topology: Let $\tau$ be the collection of all subsets $U\subset X$
that are open in the CW-topology, and  either do not contain the origin $x$,
or contain almost all subsets $\overline{e_n^2}\smallsetminus e_n^0$, which are   closed 2-cells with a 0-cell removed.
This collection of subsets obviously satisfies the axioms of a topology, and we call this topology $\tau$ the \emph{coarser topology}.
Here and throughout, almost all means all but finitely many.
The set $X$, endowed   with the coarser topology, is denoted $X_{\rm crs}$.

\begin{proposition}
\mylabel{hausdorff nonparacompact}
The space $X_{\rm crs}$ is Hausdorff but not paracompact.
\end{proposition}

\proof
Clearly, the identity map $X\ra X_{\rm crs}$ is continuous, and becomes a homeomorphism outside the
origin. Thus $X_{\rm crs}$ is Hausdorff outside the origin $x$. 
Given $y\neq x$, we choose two disjoint open neighborhoods $x\in U$, $y\in V$ on the CW-complex $X$.
By shrinking $V$, we may assume that $V$ intersects only one closed 2-cell.
By enlarging $U$, we may assume that $U$ contains all remaining closed 2-cells, while staying disjoint from $V$.
Then $U,V$ are open in the coarse topology, thus $X_{\rm crs}$ is Hausdorff.

To see that the space is not paracompact, let $U_0\subset X_\crs$ be the complement of $\bigcup_{n\geq 1} e_n^0$,
and $U_n\subset X$ by the open subset $\overline{e_n^2}\smallsetminus \left\{x\right\}$.
This gives an open covering $X_{\rm crs}=\bigcup_{n\geq 0} U_n$.
Every refinement of this covering fails to be locally finite: For each $n\geq 1$, let $U'_n$ be a
member of such a refinement that contains  $e_n^0$.
Clearly, $U'_n\subset U_n$, whence the $U'_n$ are pairwise different.
By definition of the coarser topology, each neighborhood of the origin contains almost all $\overline{e_n^2}\smallsetminus e_n^0$,
therefore intersects almost all $U'_n$.  Hence our space is not paracompact.
\qed

\begin{remark}
The space $X_{\rm crs}$ is not \emph{regular}: Consider the origin $x$ and the closed set  $A=\bigcup_{n\geq 1} e_n^0$.
Then every open neighborhood of $x$ intersects every open neighborhood of $A$.
On the other hand, the space $X_{\rm crs}$ is \emph{pointwise paracompact}, a property also called \emph{metacompactness}: 
Every open covering admits a 
refinement that is pointwise finite.  Clearly, each closed 2-cell $e_n^2\subset X_\crs$ is compact, hence
$X_\crs$ is a countable union of compacta, in other words, our space is \emph{$\sigma$-compact}.
In particular, it is  \emph{Lindel\"of}, which means that every open covering has a countable subcovering.
The reader may consult Steen and Seebach \cite{Stenn; Seebach 1978} for other counterexamples in this direction.
\end{remark}

\begin{remark}
The \emph{Kelley topology} (also called the \emph{compactly generated topology})
on a space $Y$ consists of those subsets $V\subset Y$ such that $V\cap K\subset K$
is open for each compact subset $K\subset Y$. This topology plays a  role for infinite CW-complexes, for example, to define products.
One easily checks that each compact subset $K\subset X_\crs$ is also compact with respect to the CW-topology.
From this it follows that the Kelley topology of $X_\crs$ coincides with the CW-topology.
\end{remark}

\section{\v{C}ech and sheaf cohomology}

Let $\shF$ be an abelian sheaf on a topological space $Y$. Then one has   sheaf cohomology groups $H^p(Y,\shF)$,
which are defined via global sections of  injective resolutions, and   \v{C}ech cohomology groups $\check{H}^p(Y,\shF)$, which are
computed in terms of    open coverings and local sections. The two types of cohomology groups are
related by a spectral sequence $\check{H}^p(Y,\shH^q(\shF))\Rightarrow H^{p+q}(Y,\shF)$.
For details, we refer to 
Grothendieck's exposition \cite{Grothendieck 1957} and Godement's monograph \cite{Godement 1964}.
A basic fact in sheaf theory states that the canonical map $\check{H}^1(Y,\shF)\ra H^1(Y,\shF)$
is bijective, and we have a short exact sequence
$$
0\lra \check{H}^2(Y,\shF)\lra H^2(Y,\shF)\lra \check{H}^1(Y,\shH^1(\shF))\lra 0,
$$
compare \cite{Grothendieck 1957}, page 177.
Thus \v{C}ech cohomology does not coincide with sheaf cohomology provided   $\check{H}^1(Y,\shH^1(\shF))\neq 0$.

Our task is therefore to find such a situation.
Consider the CW-complex $X$ and the space $X_{\rm crs}$ constructed in the preceding section.
The following fact will be useful:

\begin{lemma}
\mylabel{cohomology coincides}
For each open subset $V\subset X_{\rm crs}$, the sheaf cohomology groups $H^p(V,\ZZ)$, $p\geq 0$
are the same, whether computed in the CW-topology or in the coarser topology.
\end{lemma}

\proof
Let $i:X\ra X_\crs$ be the identity map, which is continuous.
We have a canonical map $\ZZ_{X_\crs}\ra i_*(\ZZ_X)$
of abelian sheaves, where the left hand side is the sheaf of locally constant integer-valued functions on $X_\crs$,
and the right hand side is the direct image sheaf of the corresponding sheaf on $X$.
We first check that this map is bijective. The question is local on $X_\crs$, and bijectivity  is obvious outside the origin.
Injectivity holds because the mapping $i$ is surjective.
Since there are arbitrarily small open neighborhoods $x\in V\subset X_\crs$ that
are pathconnected and hence connected in the CW-topology, the canonical map  $\ZZ_{X_\crs,x}\ra i_*(\ZZ_X)_x$ is bijective 
as well.

In light of the Leray--Serre spectral sequence
$$
H^p(X_\crs,R^qi_*(\ZZ_X))\Longrightarrow H^{p+q}(X,\ZZ_X),
$$
it suffices to check that  $R^pi_*(\ZZ_X)=0$ for all $p>0$. This is again local, and holds for trivial reasons
outside the origin. Since there are arbitrarily small open neighborhoods $x\in V\subset X_\crs$
that are contractible in the CW-topology, and singular cohomology coincides with 
sheaf cohomology for CW-complexes (\cite{Bredon 1967}, Chapter III, Section 1),
vanishing holds at the origin as well.
\qed

\medskip
Let $U\subset X$ be the complement of the 1-skeleton $X^1\subset X$, that is, the union of all 2-cells, and
$\ZZ_U$ be the abelian sheaf of locally constant integer-valued functions  on $U$.
Clearly, $U$ is open in the coarser topology.
Thus the inclusion map $i:U\ra X_{\rm crs}$ is continuous.
From this we obtain an abelian  sheaf $\shF=i_!(\ZZ_U)$ on $X_{\rm crs}$, called \emph{extension by zero}.
It is defined by the rule
$$
\Gamma(V,\shF)=\begin{cases}
\Gamma(V,\ZZ_U) & \text{if $V\subset U$;}\\
0               & \text{else},
\end{cases}
$$
compare \cite{SGA 2}, Expose I. Its first cohomology is easily computed:

\begin{proposition}
\mylabel{first cohomology}
Let $V\subset X_\crs$ be an open subset with $H^1(V,\ZZ)=0$.
Then we have a canonical identification
$$
H^1(V,\shF) = H^0(V\cap X^1,\ZZ)/H^0(V,\ZZ).
$$
\end{proposition}

\proof
The short exact sequence $0\ra\shF\ra\ZZ_{X_\crs}\ra\ZZ_{X^1}\ra 0$ induces a long exact sequence
$$
H^0(V,\ZZ)\lra H^0(V\cap X^1,\ZZ)\lra H^1(V,\shF)\lra H^1(V,\ZZ),
$$
and the result follows.
\qed

\medskip
For this sheaf, \v{C}ech cohomology does not coincide with sheaf cohomology, in a rather drastic way:

\begin{theorem}
\mylabel{cohomology uncountable}
For the abelian sheaf $\shF=i_!(\ZZ_U) $ on the topological space $X_{\rm crs}$ the group $\check{H}^1(X_{\rm crs},\shH^1(\shF))$
is uncountable. In particular, the inclusion $\check{H}^2(X_{\rm crs},\shF)\subset H^2(X_{\rm crs},\shF)$ is not bijective.
\end{theorem}

\proof
Let $\foU=(U_\alpha)_{\alpha\in I}$ be an open covering of $X_{\rm crs}$. By definition,
the corresponding group $\check{H}^1(\foU,\shH^1(\shF))$ 
is the first cohomology of the  complex
\begin{equation}
\label{complex}
\prod_\alpha H^1(U_\alpha,\shF) \lra 
\prod_{\alpha<\beta} H^1(U_{\alpha\beta},\shF)\lra 
\prod_{\alpha<\beta<\gamma} H^1(U_{\alpha\beta\gamma},\shF).
\end{equation}
Here we employ the usual abbreviation $U_{\alpha\beta}=U_\alpha\cap U_\beta$ et cetera.
The coboundary maps are the usual one, for example \ $(s_\alpha)\mapsto(s_\beta|U_{\alpha\beta}-s_\alpha|U_{\alpha\beta})$,
and we have chosen a total order on the index set $I$. By definition, \v{C}ech cohomology equals
$$
\check{H}^1(X,\shH^1(\shF)) = \dirlim_{\foU} \check{H}^1(\foU,\shH^1(\shF)),
$$
where the direct limit runs over all open coverings ordered by the refinement relation.
For a precise definition of the  maps in the direct system, and their well-definedness,
we refer to  \cite{Godement 1964}, Chapter II, Section 5.7.

In general, it can be difficult to control such direct limits. However, one may restrict to open
coverings forming a cofinal subsystem. Therefore, we may assume that our open covering
satisfies the following five additional assumptions:
(i) Each $U_\alpha$ and the intersection $U_\alpha\cap X^1$ are, if nonempty,  contractible in the CW-topology.
(ii) Each 0-cell is contained in precisely one $U_\alpha$.
(iii) If some $U_\alpha$ contains a 0-cell $e_n^0$, $n\geq 1$, then it is contained
in the corresponding closed 2-cell $\overline{e_n^2}$.
(iv) We  suppose that the index set $I$   is   well-ordered. This allows us  
to  regard the natural numbers $0,1,\ldots\in I$ as indices.
After reindexing, we   stipulate that $x\in U_0 $ and $e_n^0\subset U_n$.
(v) Finally, if a closed 2-cell $\overline{e_n^2}$ is contained in $U_0\cup U_n$,
then it is disjoint from all other  $U_\alpha$.

From now on, we only consider open coverings $\foU$ satisfying these five condition.
Choose $m\geq 1$ so that $U_0$ contains all  $\overline{e_n^2}\smallsetminus e_n^0$, $n\geq m$.
Condition (i) implies  that for $V=U_0\cap U_n=U_n\smallsetminus e_n^0$, $n\geq m$ we have $H^1(V,\ZZ)=0$, and furthermore
$$
H^0(V\cap X^1)=\ZZ^{\oplus 2}\quadand H^0(V,\ZZ)=\ZZ,
$$
the latter sitting diagonally in the former. Note that this is the key step in Grothendieck's argument \cite{Grothendieck 1957},
page 178. Now Proposition \ref{first cohomology} gives us an identification
$$
\prod_{n=m}^\infty H^1(U_0\cap U_n) =\prod_{n=m}^\infty\ZZ.
$$
In light of Proposition \ref{first cohomology},  Condition (i) ensures that the term on the left in the
complex (\ref{complex}) vanishes. Condition (ii) and (v) tell  us that the triple intersections $U_0\cap U_n\cap U_\alpha$ are empty
for $n\geq m$ and all indices $\alpha\neq 0,n$.
The upshot is that we have a canonical inclusion
$$
\prod_{n=m}^\infty\ZZ = \prod_{n=m}^\infty H^1(U_0\cap U_n,\shF)\subset \check{H}^1(\foU,\shH^1(\shF)).
$$
If $\foU'$ is a refinement of $\foU$ satisfying the same five conditions,
the induced map $\check{H}^1(\foU,\shH^1(\shF))\ra \check{H}^1(\foU',\shH^1(\shF))$
restricts to  the canonical projection 
$$
\prod_{n=m}^\infty \ZZ\lra\prod_{n=m'}^\infty\ZZ
$$
on the subgroups considered above, where we tacitly choose $m'\geq m$.
Since forming direct limits is exact, we obtain an inclusion
$$
\dirlim_m\prod_{n=m}^\infty\ZZ \subset \dirlim_{\foU}\check{H}^1(\foU,\shH^1(\shF))= \check{H}^1(X,\shH^1(\shF)).
$$
Again using that forming direct limits is  exact, we may rewrite the left hand side as
$$
\dirlim_{m} \left(\prod_{n=1}^\infty\ZZ \left/ \prod_{n=1}^{m-1}\ZZ\right) \right.=
\left(\prod_{n=1}^\infty\ZZ\right) \left/ \left( \dirlim_{m}\prod_{n=1}^{m-1}\ZZ \right) \right.=
\left(\prod_{n=1}^\infty\ZZ\right) \left/ \left(\bigoplus_{n=1}^\infty\ZZ\right),\right.
$$
which is uncountable.
\qed

\section{Continuous functions and principal bundles}

We finally examine  pathological properties of  continuous functions and principal bundles on $X_\crs$.
Let us write $\shC_{X_\crs}$ for the sheaf of continuous real-valued functions on $X_\crs$.

\begin{proposition}
\mylabel{nontrivial cohomology}
We have $H^1(X_\crs,\shC_{X_\crs})\neq 0$.
\end{proposition}

\proof
Recall that \v{C}ech cohomology agrees with sheaf cohomology in degree one.
Thus our task is   to construct a nontrivial \v{C}ech cohomology class.
Consider the open covering $\foU$ given by  $U_0=X_\crs\smallsetminus\bigcup_{n\geq 1}e_n^0$
and $U_n=\overline{e_n^2}\smallsetminus\left\{x\right\}$, $n\geq 1$.
Choose germs of continuous functions $f_n:(U_n,e_n^0)\ra\RR$
having an isolated zero at $e_n^0$. Then its reciprocal $1/f_n$ is defined
on some open punctured neighborhood of $e_n^0\subset U_n$, where it is necessarily  unbounded.
On the other hand, for any continuous function $g:U_0\ra\RR$  there is some $m\geq 0$ so that 
$g$ is  bounded on $\bigcup_{n=m}^\infty \overline{e_n^2}\cap U_0$.
Whence $1/f_n$ cannot be written as the difference of continuous functions
coming from $U_0$ and $U_n$, for $n\geq m$. The same applies  for any refinement $\foU'$ satisfying the five conditions
formulated in the proof for Theorem \ref{cohomology uncountable}. The upshot is that for all refinements $\foU'$ with
$U_n'$ sufficiently small, we obtain a well-defined tuple
$$
(1/f_n)_{n\geq m}\in \prod_{n=m}^\infty H^0(U'_0\cap U'_n,\shC_{X_\crs})
$$
that is a cocycle whose class in $\check{H}^1(\foU',\shC_{X_\crs})$ is nonzero. Recall that $m\geq 1$ is
any integer so that $U_0'$  contains $\overline{e_n^2}\smallsetminus e_n^0$ for all $n\geq m$.
Since this holds for all such refinements $\foU'$, 
it follows that the class  in the direct limit $\check{H}^1(X_\crs,\shC_{X_\crs})$ is nonzero as well.
\qed

\begin{remark}
On   normal spaces $Y$, the Uryson Lemma ensures that the sheaf $\shC_Y$ is \emph{soft}, that is,
the canonical map $H^0(Y,\shC_Y)\ra H^0(A,i^{-1}(\shC_Y))$ is surjective for all closed subsets $A$, where $i:A\ra Y$
denotes the inclusion map. 
According to \cite{Godement 1964}, Chapter II, Theorem 4.4.3, soft sheaves on paracompact spaces $Y$ are acyclic.

For our space $X_\crs$, it is easy to check that the canonical map 
for the discrete closed subset $A=\bigcup_{n\geq 1}e_n^0$ is not surjective.
Summing up,  the sheaf of continuous real-valued functions on $X_\crs$ is neither soft nor acyclic.
\end{remark}

Next consider the sheaf $\shS^1_{X_\crs}$   of continuous functions taking values in the circle group $S^1=\RR/\ZZ$.
It sits in the exponential sequence
$$
0\lra\ZZ\lra \shC_{X_\crs}\lra \shS^1_{X_\crs}\lra 0,
$$
where the map on the right is $t\mapsto e^{2\pi it}$. From this we get an exact sequence
$$
H^p(X_\crs,\ZZ)\lra H^p(X_\crs,\shC_{X_\crs})\lra H^p(X_\crs,\shS^1_{X_\crs})\lra H^{p+1}(X_\crs,\ZZ).
$$
The outer terms vanish by Lemma \ref{cohomology coincides}, and we conclude that 
$$
H^p(X_\crs,\shC_{X_\crs})= H^p(X_\crs,\shS^1_{X_\crs})
$$
for all $p>0$. From the preceding Proposition we get $H^1(X_\crs,\shS^1_{X_\crs})\neq 0$. In other words, there
are nontrivial principal $S^1$-bundles over $X_\crs$. This is in stark contrast to the following fact:

\begin{proposition}
\mylabel{homotopic constant}
The space $X_\crs$ is contractible.
\end{proposition}

\proof
The CW-topology and the coarse topology induce the same topology on the compact subsets
$\overline{e_n^2}\subset X_\crs$, which are thus homeomorphic to the 2-disk.
Choose homotopies $h_n:\overline{e_n^2}\times I\ra \overline{e_n^2}$ between the identity
and the constant map to the origin so that $h_n(x,t)=x$ for 
all $t\in I$, and $h_n(y,t)\not\in e_n^0$ for all $t>0$ and all $y$.
The first condition ensures that the homotopies glue to a map $h:X\times I\ra X$, which is 
continuous with respect to the CW-topology. 
From the second condition one easily infers  that it remains  continuous when regarded as 
a map $h:X_\crs\times I\ra X_\crs$. Thus $h$ is a homotopy from
the identity on $X_\crs$ to the constant map $X_\crs\ra\left\{x\right\}$.
\qed

\medskip
Let $G$ be a topological group. Recall that  a $G$-principal bundle $P\ra Y$ is called \emph{numerable} 
if it can be trivialized on some numerable covering $Y=\bigcup_{\alpha\in I} V_\alpha$.
The latter means that there is a  partition of unity $f_\beta:X\ra [0,1]$, $\beta\in J$ so that the open
covering $f^{-1}(]0,1])$, $\beta\in J$ is locally finite and refines the given covering $V_\alpha$, $\alpha\in I$.
According to Milnor's construction of the classifying space 
$$
BG=G\star G\star G\star\ldots
$$
as a countable join \cite{Milnor 1956}, together with  Dold's analysis (\cite{Dold 1963}, Section 7 and 8), 
the isomorphism classes of numerable bundles correspond to the homotopy classes of continuous maps $Y\ra BG$.
We conclude:

\begin{corollary}
\mylabel{nonnumerable}
The only principal $G$-bundles over $X_\crs$ that are numerable are the trivial ones.
\end{corollary}

\begin{remark}
Non-numerable principal $\ZZ$-bundles based on a construction with the long line  appear  in \cite{Bredon 1968}.
A non-numerable principal $\RR$-bundle over a non-Hausdorff space is sketched in  \cite{tom Dieck 2008}, page 350.
\end{remark} 

\begin{remark}
The results in this section hold true if one uses a simpler space, obtained by attaching
only 1-cells $e_n^1$ and no 2-cells, rather than pairs of 1-cells $e_{\pm n}^1$ and 2-cells $e_n^2$.
Of course, the coarser topology is defined in the same way.
\end{remark}


\end{document}